\documentclass{amsart}

\usepackage{amsfonts}
\usepackage{amsmath}
\usepackage{amsthm}
\usepackage{amssymb}
\usepackage[english]{babel}
\usepackage{graphics}
\usepackage{latexsym}
\usepackage{longtable} \usepackage{mathrsfs}
\usepackage{graphicx}
\usepackage{wrapfig} 
\usepackage[pdftex,colorlinks=true]{hyperref}

\newtheorem{theorem}{Theorem}
\newtheorem{corollary}[theorem]{Corollary}
\newtheorem{lemma}[theorem]{Lemma}

\newtheorem{claim}[theorem]{Claim}
\newtheorem{example}[theorem]{Example}

\theoremstyle{definition}
\newtheorem{definition}[theorem]{Definition}
\newtheorem{remark}[theorem]{Remark}

\newcommand{\msL}{\mathscr{L}}

\newcommand{\mC}{\mathcal{C}}

\newcommand{\mF}{\mathcal{F}}

\newcommand{\mP}{\mathcal{P}}

\newcommand{\R}{\mathbb{R}}
\newcommand{\N}{\mathbb{N}}

\newcommand{\noi}{\noindent}
\newcommand{\ms}{\medskip}

\newcommand{\al}{\alpha}

\newcommand{\e}{\varepsilon}

\newcommand{\Om}{\Omega}

\newcommand{\lharpoonup}{-\!\!\!-\!\!\!\!\rightharpoonup}
\newcommand{\larrow}{\longrightarrow}

\newcommand{\ri}{\rightarrow}

\newcommand{\p}{\partial}
\newcommand{\sub}{\subseteq}
\newcommand{\set}{\setminus}
\newcommand{\by}{\times}

\newcommand{\ess}{\textrm{ess}}
\newcommand{\diam}{\textrm{diam}}
\newcommand{\dist}{\textrm{dist}}

\newcommand{\co}{\overline{\textrm{co}}}
\newcommand{\inter}{\textrm{int}}
\newcommand{\spn}{\textrm{span}}

\renewcommand{\^}[1]{^{#1}_{\phantom{l}}}

\newcommand{\bt}{\begin{theorem}}\newcommand{\et}{\end{theorem}}
\newcommand{\bd}{\begin{definition}}\newcommand{\ed}{\end{definition}}
\newcommand{\bl}{\begin{lemma}}\newcommand{\el}{\end{lemma}}
\newcommand{\beq}{\begin{equation}}\newcommand{\eeq}{\end{equation}}
\newcommand{\bc}{\begin{claim}}\newcommand{\ec}{\end{claim}}
\newcommand{\bp}{\begin{proof}}\newcommand{\ep}{\end{proof}}

\newcommand{\BPC}{\medskip \noindent \textbf{Proof of Claim} }

\newcommand{\BPT}{\medskip \noindent \textbf{Proof of Theorem} }

\numberwithin{equation}{section} \numberwithin{theorem}{section}

\begin{document}

\title{Maximum Principles for Vectorial Approximate Minimizers of Nonconvex
Functionals \ \ }

\author{\textsl{Nikolaos I. Katzourakis}}
\address{BCAM - Basque Center for Applied Mathematics, Building 500, Biskaia Technology Park, Derio, E-48160, Spain}
\email{nkatzourakis@bcamath.org}


\date{}



\begin{abstract} We establish Maximum Principles which apply
to vectorial approximate minimizers of the general integral
functional of Calculus of Variations. Our main result is a version
of the Convex Hull Property.

The primary advance compared to results already existing in the
literature is that we have dropped the quasiconvexity assumption of
the integrand in the gradient term. The lack of weak Lower
Semicontinuity is compensated by introducing a nonlinear convergence
technique, based on the approximation of the projection onto a
convex set by reflections and on the invariance of the integrand
 in the gradient term under the Orthogonal Group.

Maximum Principles are implied for the relaxed solution in the case
of non-existence of minimizers and for minimizing solutions of the
Euler-Lagrange system of PDE.
\end{abstract}

\maketitle

\section{Introduction.} \label{Introduction}

Let $\msL : \Om \by \R^N \by \R^{N \by n} \larrow \R$ be a
Carath\'eodory function where $\Om \sub \R^n$ is an open set. In
this paper we are concerned with the derivation of Maximum Principle
results applying to \emph{Approximate Minimizers} of the functional
 \beq \label{equ1}
E(u,\Om)\ := \ \int_\Om \msL \big(x,u(x),Du(x)\big)dx
 \eeq
placed in $[W^{1,q}_g(\Om)]^N$, $q\geq 1$ with prescribed boundary
values $g\in [W^{1,q}(\Om)]^N$. A vector function $u :\Om \sub \R^n
\larrow \R^N$ will be called an \emph{$\al$-Minimizer} of
\eqref{equ1} if for some $\al\geq 0$ and all $\phi \in
[W^{1,q}_0(\Om)]^N$, we have
 \beq  \label{equ2}
E(u,\Om) \ \leq \ E(u+\phi,\Om)\ + \ \al.
 \eeq
Minimizing families $\{u_\al \}_{\al >0}$ correspond to minimizing
sequences as $\al \ri 0$ of the variational problem
 \beq \label{equ3}
E(u,\Om)\ = \  \inf_{[W^{1,q}_g(\Om)]^N}E
 \eeq
and $0$-minimizers correspond to solutions of \eqref{equ3}. Our
viewpoint of the Maximum Principle in the vector case is based on
the observation that the scalar inequalities $\sup_\Om u \leq
\max_{\p \Om} u$, $\inf_\Om u \geq \min_{\p \Om} u$ when $N=1$ can
be recast as $u(\Om) \sub \big[\min_{\p \Om} u,\ \max_{\p \Om}
u\big] $. When $N\geq1$, the appropriate vectorial extension is the
so-called \emph{Convex Hull Property}
\begin{equation}
\label{CHP}
 u(\Om)\ \sub \ \co\, \big(u(\p \Om)\big)
\end{equation}
and states that \emph{the range is contained in the closed convex
hull of the boundary values}.

Maximum Principles either as norm bounds or in the form \eqref{CHP}
applying to the functional \eqref{equ1} and its minimizing solutions
of the respective Euler-Lagrange system of PDEs are a
well-established subject and provide a priori localization and
bounds, necessary for further regularity investigations. Without
attempting to provide a complete list of papers, we refer to
Alexander-Ghomi \cite{AG}, Bildhauer-Fuchs \cite{BF1}, \cite{BF2},
Colding-Minicozzi II \cite{CM}, D'Ottavio-Leonetti-Musciano
\cite{DLM}, Leo-netti \cite{L1}, \cite{L2}, Leonetti-Siepe \cite{LS}
and Osserman \cite{O1}, \cite{O2}.

In this paper we derive appropriate versions of the Convex Hull
Property \eqref{CHP} for approximate minimizers of the functional
\eqref{equ1}. These results are inherited by the relaxed solutions
of \eqref{equ3} and by minimizing solutions of the Euler-Lagrange
system of PDE, in case they exist.

Our primary advance compared to similar results existing in the
literature is that we have dropped the quasiconvexity assumption of
$\msL$ in the gradient term. This is achieved by introducing a
nonlinear approximation technique which establishes \eqref{CHP}
without invoking the weak Lower Semi-Continuity of \eqref{equ1}. The
essential ingredients are the \emph{construction of systematic
approximations of the projection map onto convex sets generated by
reflections} and some (almost) \emph{invariance of $P \mapsto
\msL(x,\eta,P)$ under the Orthogonal Group of $\R^N$}. The latter is
automatically satisfied if $\msL=\msL(x,\eta,|P|)$.

Our approximations of the projection map by reflections, the
``Folding Maps'', are employed to construct suitable energy
comparison functions. Roughly, by projecting the minimizer on the
convex hull of its boundary values we would get an energy-decreasing
deformation. However, in the absence of convexity by projecting we
may \emph{not} decrease the energy. We obtain our deformations as
energy-preserving piecewise isometries which converge to the
projection.

A subtle point of the nonconvex case is that we can assert a version
of \eqref{CHP} only for at least one approximate minimizer among all
that may exist. This allows to utilize \eqref{CHP} as a Selecting
Principle in the sense of Dacorogna-Ferriero \cite{DF} to rule out
non-physical solutions. Moreover, our approximation method is fairly
general and could be useful in other contexts as well. For this
reason, Section \ref{Construction of the Approximations.} which is
devoted to the construction of the approximations is independent of
Section \ref{The Folding Principles.} in which the approximations
are employed to derive our Maximum Principles. We also note that although 
our deformations are ``non-physical''  since they are constructed via reflections, our
results could be useful to nonconvex problems arising
in Elasticity since rotations can be factored as compositions of reflections.

The main ideas arising in this work is an outgrowth of methods
employed as stepping stones in order to solve the problem studied in
\cite{AK}, by Alikakos and the Author. The choice to formulate our
results for approximate minimizers is essential; in the absence of
quasiconvexity, existence of minimizers is not to be expected.
Moreover, approximate minimizers satisfy a locality property which
is not enjoyed by similar notions of ``almost'' minimizers, like
$Q$-minimizers or $\omega$-minimizers (cf.\ Giusti \cite{Gi},
Dacorogna \cite{D1}). This property is crucial for our methods and
owes to the fact that deviation form minimality is viewed
additively.

This paper is organized as follows. In Section \ref{Construction of
the Approximations.} we present the core of our method, which
consists of the construction of the approximations. Then, in Section
\ref{The Folding Principles.} we employ them to establish our main
results, the Maximum Principles.

\section{Construction of the Approximations.} \label{Construction of the Approximations.}

Let us begin with some basics. In what follows, Sobolev functions
will always be tacitly identified with their precise
representatives. The Orthogonal group of $\R^N$ will be denoted by
$O(N,\R)$ and the \emph{Affine} Orthogonal group of $\R^N$ by
$AO(N,\R)$.

Let $\xi \in \mathbb{S}^{N-1}$ be a vector in the unit sphere. $\xi$
defines a hyperplane $H_{\xi}:=(\spn [\xi])^\bot \subseteq \R^N$ and
a reflection $R_{\xi} \in \textrm{O(N},\R)$ with respect to
$H_{\xi}$, given by $R_{\xi}u:=u-2(u\cdot\xi)\xi$. Here, ``$\cdot$''
denotes inner product. If $u_0 \in \R^N$, the affine reflection with
respect to the hyperplane $H_{\xi}+u_0$ is given by $u \mapsto u-
2\big((u-u_0)\cdot \xi)\xi$ $= R_\xi u+2(u_0\cdot \xi)\xi$.

Let $\mC \sub \R^N$ be compact and convex. We will write
$\xi=\hat{u}$ for the outwards pointing normal to the supporting
hyperplane $H_{\hat{u}} +u$ at $u$, but unless $\p \mC$ is
differentiable, $u\mapsto \hat{u}$ may be multivalued and determines
the normal cone. For any $u \in \p\mC$, $\hat{u}$ generates a
reflection $w \mapsto R_{\hat{u}}w +2(u \cdot \hat{u})\hat{u}$ in
$\textrm{AO}(\textrm{N},\R)$ with respect to $H_{\hat{u}}+u$.

The following is the principal result of this section. It states
that given any compact convex set, there exists a sequence of
weak$^*$ approximations in $W^{1,\infty}_{loc}(\R^N)^N$ of the
projection on it generated by affine reflections \emph{merely}.

\begin{theorem}[Approximation of the projection onto convex sets by reflections] \label{Construction of the
Approximations}

\noi Let $\mC \sub \R^N$ be a compact convex set with nonempty
interior and $0 \in \textrm{int}(\mC)$. Let also $\mP^\mC:\R^N
\larrow \mC$ denote the projection on it.

Then, there exists a sequence of locally Lipschitz maps $\mF^\mC_m
:\R^N \larrow \R^N$, $m\in \N$, which satisfies:

\noi(i) Each $\mF^\mC_m$ is piecewise equal to a composition of
finitely many successive affine reflections. In particular, for
a.e.\ $u\in \R^N$, there is an $R \in O(N,\R)$ such that $D
\mF^\mC_m(u)=R$.

\noi (ii) Each $\mF^\mC_m$ equals the identity on $\mC$ and maps the
$m$-dilate of $\mC$ into the $(1+\frac{1}{m})$-dilate of $\mC$:
 \beq \label{eq1}
\left\{
 \begin{array}{l}
\mF^\mC_m (E)\ = \  E, \ \ E\sub \mC, \medskip\\
\mF^\mC_m(m \mC) \ \sub \ \left(1+\dfrac{1}{m}\right)\mC.
\end{array}
\right.
 \eeq

\noi(iii) The sequence $\mF^\mC_m$ is a locally uniform
approximation of the projection onto the convex set $\mC$; moreover,
we have $\mF^\mC_m \ \overset{*}{\lharpoonup} \ \mP^\mC$ weakly$^*$
in $[W^{1,\infty}_{\textrm{loc}} (\R^N)]^N$, as $m \rightarrow
\infty$.
\end{theorem}

\noi The functions $\{\mF_m^\mC\ |\ m \in \N\}$ provided by Theorem
\ref{Construction of the Approximations} will be referred to as the
\emph{Folding Maps of the convex set $\mC$}.

\begin{example} An elementary example of Folding maps in the scalar
case is given by the sequence defined inductively as $\mF_0(u):=u
\chi_{(0,1)}(u)+(2-u) \chi_{(1,2)}(u)$ and $\mF_{m+1}:= \mF_m
\chi_{\{\mF_m < 2\^{-(m+1)}\}} +(2\^{-m}-\mF_m)\chi_{\{\mF_m >
2\^{-(m+1)}\}} $, in $W\^{1,\infty}(0,2)$.
\[
\includegraphics[scale=0.17]{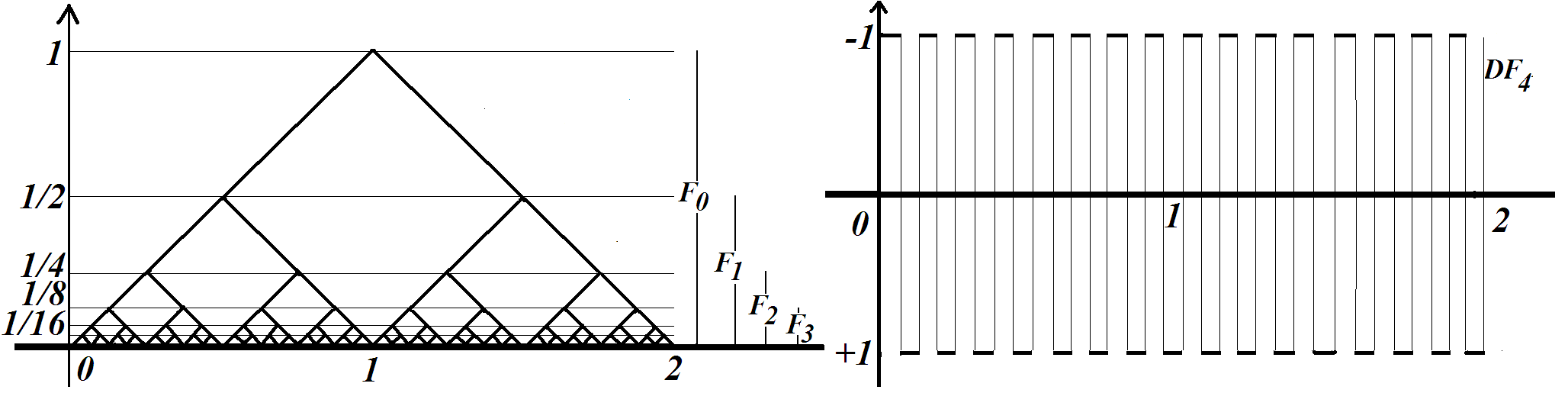}
\]

We have that $ \|\mF_m\|_{W\^{1,\infty} (0,2)}\leq 1$ and since
$\|\mF_m\|_{L\^{\infty}(0,2)} \leq 2\^{-m-1}$, the weak$^*$ limit is
zero. Thus $\mF_m \overset{*}{\lharpoonup} 0$ in
$W\^{1,\infty}(0,2)$, which implies $\mF_m \larrow 0$ in
$L\^{\infty}(0,2)$. The sequence $\mF_m$ approximates the projection
$\mP^0$ on $\{0\}$, which is the zero map.
\end{example}

\noi Throughout the proof we will need to supplement convexity with
an auxiliary notion which rectifies the lack of smoothness. It is
eventually dropped by approximation.

Let $\mathcal{C}$ be a \emph{simplex}, that is a compact set of the
form $\bigcap \big\{H^-_j : 1 \leq j \leq K\big\}$ generated by
closed halfspaces $H^-_j$ which are determined by affine hyperplanes
$H_j=(\spn [\xi_i])^{\bot}+u_i$. For each hyperplane $H_i$, $\xi_i$
is its unit normal vector pointing outwards of the halfspace
$H^-_i$, while the hyperplane passes through $u_i$.

A simplex $\mathcal{C} \sub \R^N$ will be called \textbf{obtuse}, if
its defining sides $H^{*}_j$ meet in obtuse angles $\geq
\dfrac{\pi}{2}$, that is, when for all $i, j \in \{1,...,K\}$ we
have the implication
\begin{equation} \label{O}
\ H^{*}_i \cap H^{*}_j \ \neq \ \emptyset \ \ \ \Rightarrow\ \ \ 0
\leq \xi_i \cdot \xi_j \leq 1.
\end{equation}

\BPT \ref{Construction of the Approximations}. Since $\mC$ is fixed,
we drop the superscript ``$\mC$''. We shall first prove the
existence of $\mF_m$ for obtuse simplices satisfying \eqref{O} and
then deduce the case of general convex sets by approximation. For,
assume $\mC =\bigcap \{ H^-_j : 1 \leq j \leq K \}$ with $K\geq N+1$
and $H_j=u_j +(\spn [\xi_j])^\bot$ with \eqref{O} being satisfied by
the sides $H^*_i$. We fix $m \in \N$. Then, $\mC \sub
\left(1+\frac{1}{m}\right) \mC$ and since
$\left(1+\frac{1}{m}\right) H_j=\left(1+\frac{1}{m}\right) u_j +
(\spn [\xi_j])^\bot$, we obtain
 \beq
\left(1+\frac{1}{m}\right) \mC\ =\ \bigcap
\left\{\left(1+\frac{1}{m}\right)  H^-_j : 1 \leq j \leq K \right\}.
 \eeq
We choose $\xi_{\min}$ on $\p \mC$ and $\xi_{\max}$ on $\bigcup \{
H_j : 1 \leq j \leq K \}$ such that
\begin{align}
|\xi_{\min}|\  &= \ \min_{1 \leq j \leq K}\ \left(\min_{u \in H_j}\
|u|\right), \medskip \label{minmax1}\\
| \xi_{\max}|\  &= \ \max_{1 \leq j \leq K}\ \left(\min_{u \in H_j}\
|u|\right). \label{minmax2}
\end{align}
\[
\includegraphics[scale=0.12]{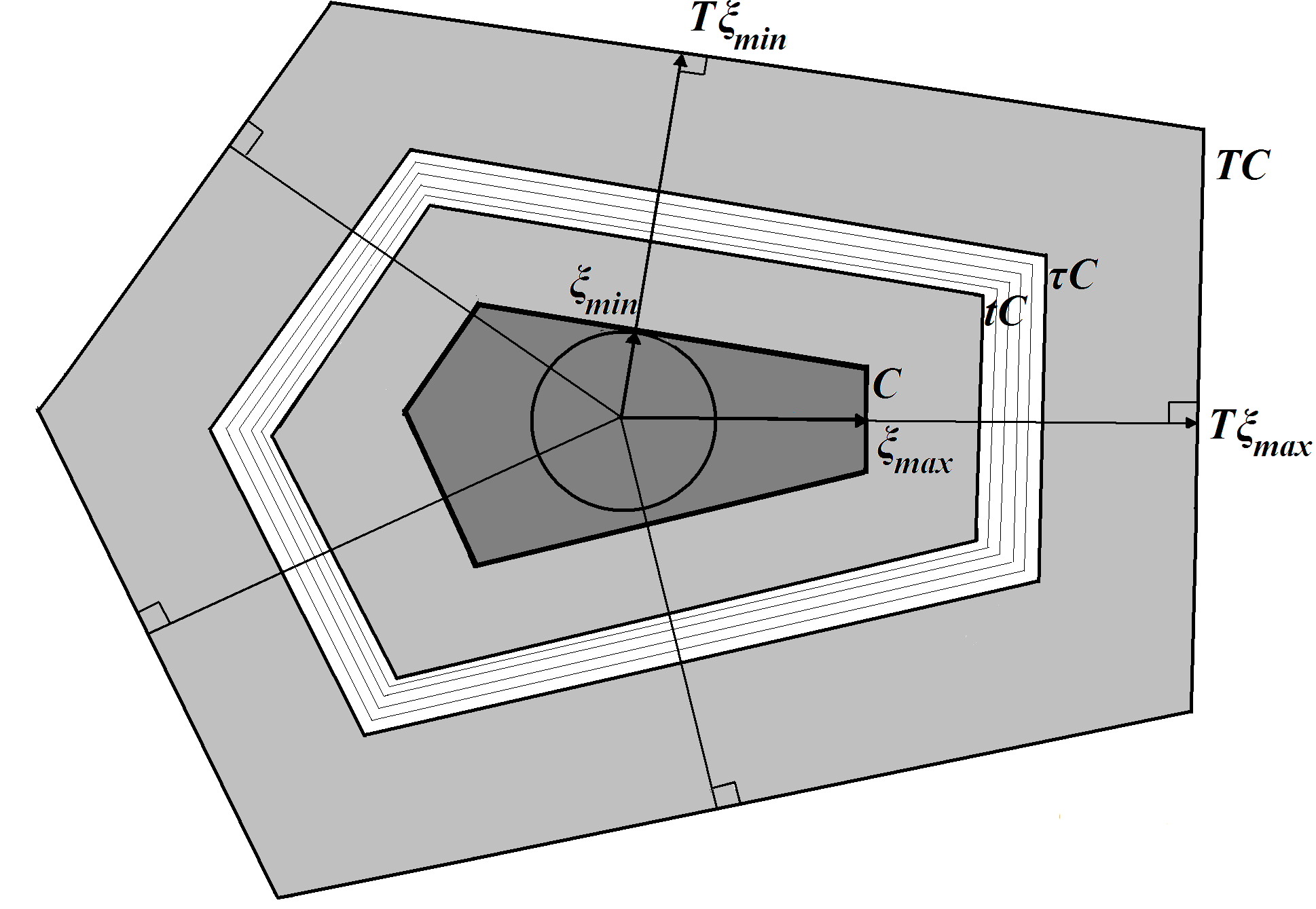}
\]
The points $\xi_{\min}$ and $\xi_{\max}$ always exist and are normal
to $\bigcup \{ H_j : 1 \leq j \leq K \}$, but possibly they are not
unique or perhaps $|\xi_{\min}|=|\xi_{\max}|$. Let us set $T:=m$.
The dilation $u\mapsto Tu$ maps $\mC$ onto $T\mC$. We fix two
numbers $t,\tau
>0$ with $1< t<\tau\leq T$, chosen such that
\begin{equation}
\label{$**$} t \ = \ \dfrac{\tau\ |\xi_{\max}|\ +\
|\xi_{\min}|}{|\xi_{\min}|\ + \ |\xi_{\max}|}.
\end{equation}
By \eqref{minmax1} and \eqref{minmax2}, the minimum distance between
the sides of $\p \mC$ and $\p (t\mC)$ is realized along the
$\xi_{\min}$ direction, while the maximum distance between the sides
of $\p (\tau \mC)$ and $\p (t\mC)$ along the $\xi_{\max}$ direction.
Thus, \eqref{$**$} rearranged says $|(t-1)\xi_{\min}| =  |(\tau
-t)\xi_{\max}|$, or
\begin{equation}
\min_{1 \leq j \leq K}\ \left(\min_{u \in H_j}\ |tu-u| \right)\ = \
\max_{1 \leq j \leq K}\ \left(\min_{u \in H_j}\ |\tau u-tu|\right).
\end{equation}
Let us fix $s \in [1,T]$ and $j \in \{1,...,K\}$. The affine
reflection $R^s_j \in \textrm{AO(N},\R)$ with respect to $sH_j$ is
given by $R^s_j(u)=R_{\xi_j}u+2s(u_j \cdot \xi_j)\xi_j$, $R_{\xi_j}
\in O(N,\R)$. For any $s \in [1,T]$ and $j \in \{1,...,K\}$, we
define
\begin{equation} \label{FI}
\mF^s_j(u)\ :=\ \left\{\begin{array}{l}
                         u\ , \ \ \ \ \ \ \ \ \ u \in sH^-_j \\
                         R^s_j(u)\ , \ \ \ u \notin sH^-_j
                       \end{array}
\right. , \ \ \ \ \ \ \ \ \   \mF_j^s \ : \ \R^N \larrow \R^N,
\end{equation}
and set
\begin{equation}\label{FII}
\mF^s \ :=\ \mF^s_K \circ \cdots \circ \mF^s_1,\ \ \ \ \ \ \ \ \mF^s
\ : \ \R^N \larrow \R^N.
\end{equation}
We obtain that $\mF^s$ equals the identity on ${s\mC}$ and for a.e.\
$u \in \R^N$, there is an $R \in O(N,\R)$ such that $D\mF^s (u)=Ru$.
Hence, $|D\mF^s|=1$ a.e.\ on $\R^N$.
\begin{claim} \label{Cl1} Assume that $t,\tau$ are in $[1,T]$ and they satisfy
\eqref{$**$}. Then, for all $j=1,... ,K$, we have
\begin{equation}
\mF^t_j\Big(\tau \mC \set \inter (tH^-_j) \Big) \ \sub \ \tau \mC
\cap \Big(t H^-_j \set  \inter(H^-_j)\Big).
\end{equation}
\end{claim}

 \BPC \ref{Cl1}. First note that by \eqref{$**$},
 \beq
 \dist
(\tau H^-_j,t H^-_j) \leq \dist (t H^-_j,H^-_j),
 \eeq
for $j=1,..,K$. Thus, $\mF^t_j$ maps the zone $\tau H^-_j \set
\inter(tH^-_j)$ into the zone $t H^-_j \setminus \inter(H^-_j)$. Let
us normalize up to an element of $\textrm{AO(N},\R)$ to $\xi_j=e_N$,
$tH_j=\{u_N=0\}$. Let $v=(v',v_N)$ be in $\tau \mC \set \inter
(tH^-_j)$. Then, we have $R^t_j(v)=(v',-v_N)$ and $R^t_j(v)$ belongs
to
 \beq
tH^-_j \set \inter(H^-_j)=\{v\in \R^N: (1-t)|u_j| \leq v_N \leq 0\}.
 \eeq
By obtuseness, it suffices to show that $R^t_j(v)$ is into
 \beq
\bigcap\{\tau H^-_i : i \neq j, H^{*}_i \cap H^{*}_j \neq
\emptyset\},
 \eeq
for all the neighbor sides $H^{*}_i$ of $H^{*}_j$. Let us fix such
an $H^{*}_i$. Then, $\tau H_i=(\spn [\xi_i])^\bot+\tau u_i$ and by a
further translation we may assume $u_i = 0$ (of course, the rest
points $\{u_l : l \neq i\}$ relative to $H^*_l$ necessarily change,
but this affects nothing).
\[
\includegraphics[scale=0.12]{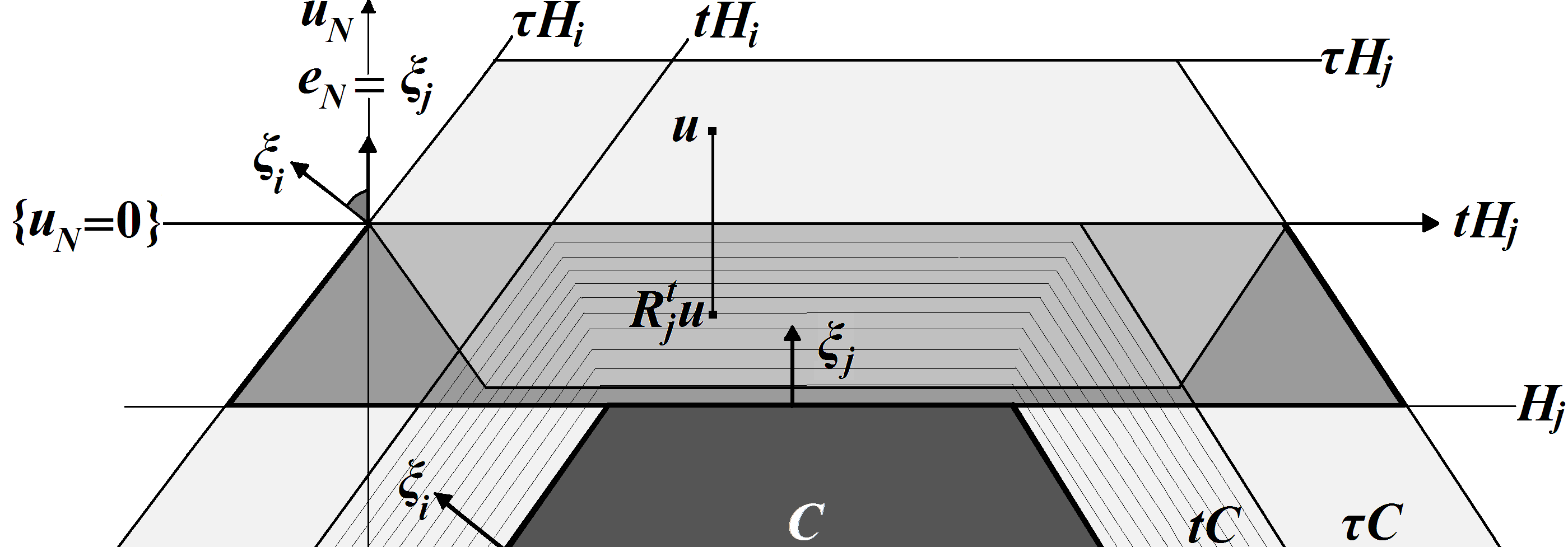}
\]
We obtain
 \beq
 \begin{array}{c}
\tau H_i = \{z \in \R^N \ :\ z \cdot \xi_i=0\}, \\
\tau H^-_i = \{z \in \R^N \ :\ z \cdot \xi_i \leq 0\}.
 \end{array}
 \eeq
Since $v \in \tau \mC$ and
 \beq
\tau \mC = \bigcap \{\tau H^-_j : 1 \leq j \leq K\} \sub \tau H^-_i,
 \eeq
we get that $v \cdot \xi_i \leq 0$ and by writing $v=(v',v_N)$ and
$\xi_i=(\xi'_i,\xi'_N)$ and denoting inner products both in $\R^N$
and $\R^{N-1}$ by ``$\cdot$'', we have
\begin{align}
u' \cdot \xi_i' \ & = \ -v'_N \xi'_{iN} \ + \ v \cdot \xi_i \nonumber\\
& = \ -(v\cdot e_N)(\xi_i \cdot e_N) \ + \ v \cdot
\xi_i\\
\nonumber & \leq \ -(v\cdot e_N)(\xi_i \cdot e_N).
\end{align}
Thus, utilizing that $v_N \geq 0$ and that obtuseness is equivalent
to $0\leq \xi_i \cdot e_N  \leq 1$ for all such $i \neq j$, we
obtain
\begin{align}
R^t_j(v) \cdot \xi_i \ & = v' \cdot \xi'_i\ + \
(-v_N)\xi'_N \nonumber\\
& = \ v' \cdot \xi_i' \ - \ (v\cdot e_N)(\xi_i \cdot e_N)\\
 & \leq \ -2(v\cdot e_N)(\xi_i \cdot e_N) \nonumber\\
 &  \leq \ 0. \nonumber
\end{align}
Consequently, if $v \in \tau \mC \set \inter (tH^-_j)$, we have
$\mF^t_j(v) \in \tau H^-_i$ for all $i \neq j$ for which $H^{*}_i
\cap H^{*}_j \neq \emptyset$. \qed

\begin{claim} \label{Cl2} Assume that $t,\tau$ belong to $[1,T]$
and that they satisfy \eqref{$**$}. Then, we have
\begin{equation}
\left\{
\begin{array}{l}
 (i)\ \ \mF^t\big(\tau \mC \set \inter(t\mC)\big)
\ \sub \ t \mC \set  \inter(\mC), \medskip\\
(ii)\ \ \mF^t \big(t \mC \big) \ = \ t \mC.
 \end{array} \right.
\end{equation}
\end{claim}
\BPC \ref{Cl2}. Let $u \in \tau \mC \set \inter(t\mC)$ be a fixed
point. We first show that $\mF^t(u) \notin \inter(\mC)$. Let $k \in
\{1,...,K\}$ be the first in order index for which $u \in \tau \mC
\set \inter(tH^-_k)$. By \eqref{FI}, all
$\mF^t_1,$,...,$\mF^t_{k-1}$ leave $u$ invariant, hence $\mF^t(u) =
\mF^t_K(...\mF^t_k(u))$. By Claim \ref{Cl1},
 \beq
\mF^t_k(u)\ \in\ \tau \mC \cap\big(tH^-_k \set \inter(H^-_k) \big)
 \eeq
and in particular $\mF^t(u) \notin \inter(\mC)$. Let $l\geq k+1$ be
the next in order index for which
\begin{equation}
\mF^t_k(u)\ \in\ \tau \mC \cap \big(tH^-_k \set \inter(H^-_k) \big)
\set \inter(tH^-_l).
\end{equation}
For all the intermediate indices $k+1 \leq i <l$, $\mF^t_i$ leaves
$\mF^t_k(u)$ invariant and hence
\begin{equation}
\mF^t_{l}\big(\mF^t_{l-1}... \mF^t_{i}...(\mF^t_{k}(u))\big)\ =\
\mF^t_{l}(\mF^t_{k}(u)).
\end{equation}
By Claim \ref{Cl1}, we have
 \beq
\mF^t_{l}(\mF^t_{k}(u))\ \in\ tH^-_l \setminus \textrm{int}(H^-_l),
 \eeq
and thus $\mF^t_{l}(\mF^t_{k}(u)) \notin \inter(\mC)$. In view on
\eqref{FII}, by continuing till $K$ we obtain $\mF^t(u) \notin
\inter(\mC)$.

\noi We now show that $\mF^t(u) \in \inter(t\mC)$. Let again $k \in
\{1,...,K\}$ be the first in order index for which $u \in \tau \mC
\setminus tH^-_k$. Then, $\mF^t(u) = \mF^t_{K}(... \mF^t_{k}(u))$.
Let $H^{*}_l$, $l\geq k+1$, be the first in order neighbor side of
$H^{*}_k$. By Claim \ref{Cl1}, we have
 \beq
 \mF^t_k(u)\ \in\ \tau \mC \cap \big(tH^-_k \set \inter(H^-_k)\big).
 \eeq
If $\mF^t_k(u) \in tH^-_l \cap \tau \mC$, we proceed to the next
neighbor side. Otherwise, if $\mF^t_k(u) \in \tau \mC \set
\inter(tH^-_l)$, the exact same argument in the proof of Claim
\ref{Cl1} with $tH_k$ in the place of $tH_j$ and $tH_l$ in the place
of $\tau H_i$ shows that
\begin{equation}
\mF^t_l(\mF^t_k(u))\ \in \ \big(tH^-_k \setminus
 \textrm{int}(H^-_k)\big) \cap \big(tH^-_l \setminus
 \textrm{int}(H^-_l)\big).
\end{equation}
In view of \eqref{FII}, continuing till $K$ we obtain
 \beq
\mF^t(u)\ \in \ \bigcap\big\{tH^-_j \set \inter(H^-_j) \ :\ 1\leq
j\leq K \big\}
 \eeq
and hence $\mF^t(u) \in t\mC$, as desired.

\noi (ii) follows directly from \eqref{FI} and \eqref{FII}.
 \qed

\ms

\noi Now we iterate Claim \ref{Cl2}. We define inductively:
\begin{equation}
\label{$***$} \left\{
\begin{array}{l}
\ \ t_0\ :=\ T\\
t_{k+1} \ := \ \dfrac{t_k \ |\xi_{\max}|\ +\
|\xi_{\min}|}{|\xi_{\min}|\ +\ |\xi_{\max}|}.
\end{array}
 \right.
\end{equation}
\begin{claim} \label{Cl3} The sequence $(t_k)^\infty_1$ defined by \eqref{$***$} is
strictly decreasing to $1^+$ as $k \rightarrow \infty$, and
\begin{equation}
t_k\ = \
\left(\dfrac{|\xi_{\max}|}{|\xi_{\max}|+|\xi_{\min}|}\right)^k T \ +
\ \left(\dfrac{|\xi_{\min}|}{|\xi_{\max}|+|\xi_{\min}|}\right)
\underset{j=0}{\overset{k-1}{\sum}}
\left(\dfrac{|\xi_{\max}|}{|\xi_{\max}|+|\xi_{\min}|}\right)^j,
\end{equation}
for all $k\in \N$.
\end{claim}

\BPC \ref{Cl3}. Follows by induction and the geometric series. \qed

\ms

\noi For $m \in \N$, we choose $k(m) \in \N$ such that $t_{k(m)} <
1+\frac{1}{m}$. Then, $t_{k(m)}\mC \sub \left(1+\frac{1}{m}\right)
\mC$. We define:
\begin{equation} \label{FIII}
\mF^\mC_m\ := \ \mF^{t_{k(m)}} \circ \cdots \circ \mF^{t_{1}}.
\end{equation}

\begin{claim} \label{Cl4} The Folding Map $\mF^\mC_m : \R^N \larrow \R^N$ given by \eqref{FI},
\eqref{FII}, \eqref{FIII} is a locally Lipschitz map in
$[W^{1,\infty}_{\text{loc}}(\R^N)]^N$, piecewise equal to
compositions of affine reflections and satisfies \eqref{eq1} for all
$m\in \N$.
\end{claim}

\BPC \ref{Cl4}. By Claims \ref{Cl2}, \ref{Cl3}, for all $k=1,2,...$
we have
\begin{equation}
\left\{\begin{array}{l} \mF^{t_k}\big(t_{k-1} \mC \set \inter(t_k \mC)\big)
 \ \sub \ t_k \mC \set  \inter(\mC), \ms\\
\mF^{t_k} \big(t_{k} \mC \big) \ = \ t_k \mC.
 \end{array} \right.
\end{equation}
By Definitions \eqref{FI}, \eqref{FII}, \eqref{FIII}, we readily
have that $\mF^\mC_m$ leaves $\mC$ invariant; thus, $\mF^\mC_m(E)=E$
for all $E\sub \mC$. The rest properties of $\mF^\mC_m$ also follow
by construction, so it suffices to establish \eqref{eq1}. For, by
employing that $m=T=t_0$, we obtain
\begin{align}
 \mF^\mC_m(m\mC) \ & \sub \
\mF^{t_{k(m)}}\big(... \mF^{t_{1}}(t_0\mC)\big) \nonumber\\
& \sub \ \mF^{t_{k(m)}}\big(... \mF^{t_{2}}(t_1\mC)\big)
\nonumber\\
& \sub  \  \mF^{t_{k(m)}}\big(... \mF^{t_{3}}(t_2\mC)\big) \\
& \sub\ ... \nonumber\\
& \sub \  t_{k(m)}\mC \nonumber\\
 & \sub \ \left(1+\frac{1}{m}\right) \mC. \nonumber
\end{align}
Now, since $T\mC \set \inter(\mC) = t_0 \mC \set \inter(\mC)$, by
utilizing \eqref{FI}, \eqref{FII}, \eqref{FIII} we have
\begin{align}
 \hspace{20pt} \mF^{t_1}(m\mC \set \mC)\ & \sub \ \mF^{t_1}\big(t_0 \mC \set
\inter(\mC)\big) \nonumber\\
& = \ \mF^{t_1} \Big(\big[(t_0 \mC \set  \inter(\mC)) \cap
\inter(t_1 \mC)
\big] \nonumber\\
& \ \ \ \ \ \cup \big[ (t_0 \mC \set  \inter(\mC)) \set \inter(t_1 \mC)\big]\Big) \\
 & \sub \ \big(t_1 \mC \set \inter(\mC)\big)\cup \mF^{t_1}\big(t_0 \mC \set \inter(t_1
 \mC)\big) \nonumber\\
 & \sub \ \big(t_1 \mC \set \inter(\mC)\big)\cup
 \big(t_1 \mC \set \inter(\mC)\big) \nonumber\\
& = \ t_1 \mC \set \inter(\mC), \nonumber
\end{align}
\begin{align}
\mF^{t_2}\big(\mF^{t_1}(m\mC \set \mC)\big)\ & \sub \
\mF^{t_2}\big(t_1
\mC \set \inter(\mC)\big) \nonumber\\
& = \ \mF^{t_2} \Big(\big[(t_1 \mC \set  \inter(\mC)) \cap \inter(t_2 \mC) \big]
\nonumber\\
& \ \ \ \ \  \cup
\big[ (t_1 \mC \set  \inter(\mC)) \set \inter(t_2 \mC)\big]\Big)\\
 & \sub \ \big(t_2 \mC \set \inter(\mC)\big)\cup \mF^{t_2}\big(t_1 \mC \set
 \inter(t_2 \mC)\big)  \nonumber\\
 & \sub \ \big(t_2 \mC \set \inter(\mC)\big)\cup
 \big(t_2 \mC \set \inter(\mC)\big) \nonumber \\
& = \ t_2 \mC \set \inter(\mC), \nonumber
\end{align}
\[
\vdots  \hspace{20pt}  \hspace{20pt}  \hspace{20pt}  \hspace{20pt}
\]
\begin{align}
\mF^{\mC}_m (m\mC \set \mC)\ & = \
\mF^{t_{k(m)}}\big(...\mF^{t_1}(m\mC \set \inter(\mC))\big)\ \ \ \ \
 \nonumber\\
& \sub \ ... \nonumber\\
& \sub \ t_{k(m)} \mC \set
\inter(\mC) \hspace{30pt} \\
 & \sub \ \left(1+\frac{1}{m}\right) \mC \set \inter(\mC) \nonumber\\
& \sub \ \left(1+\frac{1}{m}\right) \mC. \nonumber
\end{align}
Claim \ref{Cl4} has been established. \qed

\noi So far, we have established $(i)$ and $(ii)$ of Theorem
\ref{Construction of the Approximations} under the assumption that
$\mC$ is an obtuse simplex. Let us now drop this assumption.

 \bc \label{Cl5} Statements $(i)$ and $(ii)$ of Theorem
\ref{Construction of the Approximations} hold for a general compact
convex set $\mC$ with $0 \in \inter(\mC)$.
 \ec

\BPC \ref{Cl5}. Let $(\mC_k)_{k=1}^\infty$ be a decreasing sequence
of $C^1$ compact convex sets which approximates $\mC$ from the
outside:
 \beq
 \mC_1 \ \supseteq \ \mC_2 \ \supseteq \ ... \ \supseteq \ \mC_k \supseteq
  \ ... \ \supseteq \ \mC.
 \eeq
Fix $k\in\N$. Since $\mC_k$ is $C^1$, if $p,q \in \p \mC$ and
$\hat{p}$, $\hat{q}$ are the outwards pointing normal vectors at
$p,q$, then regularity and compactness imply the existence of a
continuous increasing modulus of continuity $\omega \in
C^0[0,\infty)$ with $\omega(0)=0$ such that
 \beq \label{eq}
\big| \hat{p} - \hat{q} \big| \ \leq \ \omega(|p-q|).
 \eeq
By employing estimate \eqref{eq}, we obtain
 \begin{align} \label{eq2}
\big|1\ - \ \cos(Angl(\hat{p},\hat{q}))\big|\ & = \ \left|
\frac{1}{2}\big(|\hat{p}|^2 + |\hat{q}|^2\big) - |\hat{p}|
|\hat{q}|\cos(Angl(\hat{p},\hat{q})) \right| \nonumber\\
&= \ \frac{1}{2} \left| |\hat{p}|^2 + |\hat{q}|^2 -
2(\hat{p} \cdot \hat{q}) \right|\\
&= \ \frac{1}{2} \left|
\hat{p} - \hat{q} \right|^2 \nonumber\\
& \leq \ \frac{1}{2} \omega^2(|p-q|).  \nonumber
 \end{align}
Hence, by \eqref{eq2}, if $p$ is close to $q$, the tangent
hyperplanes at these points meet at obtuse angles. By \eqref{eq2},
each smooth convex set $\mC_k$ can be approximated from the outside
by a decreasing sequence $(\mC_{k,l})_{l=1}^\infty$ of obtuse
simplices, generated by tangent hyperplanes:
 \beq
 \mC_{k,1} \ \supseteq \ \mC_{k,2} \ \supseteq \ ... \ \supseteq \ \mC_{k,l} \supseteq
  \ ... \ \supseteq \ \mC_k.
 \eeq
Without harming generality, we may assume that $l$ denotes the
number of sides of the simplex $\mC_{k,l}$. If $m \in \N$, let us
define the Folding Map $\mF_m^\mC$ of $\mC$ as the $l(m)$-Folding
Map referring to the obtuse simplex $\mC_{k(m),l(m)}$ for some
$k(m)$, $l(m) \in \N$ sufficiently large:
 \beq \label{FIV}
\mF^\mC_m \ :=\ \mF_{l(m)}^{\mC_{k(m),l(m)}}.
 \eeq
By the previous analysis, such numbers $k(m)$, $l(m)$ exist and
$\mF_m^\mC : \R^N \larrow \R^N$ has all the desired properties of
$(i)$ and $(ii)$. \qed

\ms

Let us conclude Theorem \ref{Construction of the Approximations} by
establishing $(iii)$.

 \bc \label{Cl6} There exists subsequences $(k(m))_{m=1}^\infty$,
 $(l(m))_{m=1}^\infty$ such that the Folding map of $\mC$ given by
 \eqref{FI}, \eqref{FII}, \eqref{FIII} and \eqref{FIV} satisfies
$\mF^\mC_m \ \overset{*}{\lharpoonup} \ \mP^\mC$ weakly$^*$ in
$[W^{1,\infty}_{\textrm{loc}} (\R^N)]^N$, along a subsequence as $m
\rightarrow \infty$.
 \ec

\BPC \ref{Cl6}. By Claims \ref{Cl1} - \ref{Cl5}, for any $m\in \N$
and any $j\leq m$, we have
 \begin{align}
\mF^\mC_m\big(j\mC \set \mC\big) \ & \sub \
\left(1+\frac{1}{m}\right)\mC \set
\inter(\mC), \label{eqa}\\
\mF^\mC_m(E)\ & =  \ E, \ \ E\sub \mC,  \label{eq5}\\
\big|D \mF^\mC_m\big| \ & = \ 1, \text{ a.e. on }\R^N.  \label{eq6}
 \end{align}
Hence, by \eqref{eqa} - \eqref{eq6} we obtain the bounds
 \begin{align}
\big\| \mF^\mC_m \big\|_{L^\infty(j\mC)}\ & \leq \ 2 \, \diam(\mC), \label{eq7}\\
\big\| D \mF^\mC_m \big\|_{L^\infty(\R^N)}\ & \leq \ 1,  \label{eq8}
\end{align}
for all $j\leq m$ and $m\in \N$. By \eqref{eq7}, \eqref{eq8},
weak$^*$ compactness of the local space
$[W^{1,\infty}_{\textrm{loc}} (\R^N)]^N$ and local compactness of
the imbedding $W^{1,\infty}_{\textrm{loc}} (\R^N) \subset \subset
L^\infty_{\textrm{loc}} (\R^N)$, we can extract a subsequence
(denoted again by $\mF^\mC_m$) such that, for some $\mF^\mC$ in
$[W^{1,\infty}_{\textrm{loc}} (\R^N)]^N$
 \begin{align}
 \mF^\mC_m \ & \larrow \ \mF^\mC, \text{ in }
[L^{\infty}_{\textrm{loc}} (\R^N)]^N, \label{eq15}\\
 D\mF^\mC_m \ & \overset{*}{\lharpoonup}\ D\mF^\mC, \text{ in }
[L^{\infty}_{\textrm{loc}} (\R^N)]^{N\by N},
\end{align}
as $m\ri \infty$. By passing to the locally uniform limit to
\eqref{eqa} and \eqref{eq5} and then letting $j\ri \infty$, we
obtain
 \begin{align}
& \mF^\mC\big(\R^N \set \mC\big) \  \sub \ \p \mC, \\
& \mF^\mC(E)\  =  \ E, \ \ E\sub \mC.
 \end{align}
Let us complete Theorem \ref{Construction of the Approximations} by
establishing that for an appropriate choice of $k(m)$ and $l(m)$,
the limit $\mF^\mC$ coincides with the projection $\mP^\mC$ on
$\mC$. Fix $j,m \in \N$ and $\e>0$. Let $\mP^{\mC_{k(m)}}$ be the
projection on $\mC_{k(m)}$. Consider the estimate
\begin{align} \label{eq9}
\big\|\mF^\mC - \mP^\mC \big\|_{L^\infty (j\mC)} \ \leq & \
 \big\|\mF^\mC - \mF^\mC_m \big\|_{L^\infty (j\mC)}
 \ + \ \big\|\mF^\mC_m - \mP^{\mC_{k(m)}} \big\|_{L^\infty (j\mC)} \nonumber\\
 & + \ \big\| \mP^{\mC_{k(m)}} - \mP^\mC \big\|_{L^\infty (j\mC)}.
\end{align}
By the estimate
\begin{align} \label{eq10}
\big\| \mP^{\mC_{k(m)}}  - \mP^\mC \big\|_{L^\infty (j\mC)} \leq \
\max\Big\{ |p-p_m| & \ : p \in \p\mC,\ p_m \in \mC_{k(m)},\\
& [p,p_m ] \sub \mC_{k(m)}\set \inter(\mC) \Big\},\nonumber
\end{align}
if $k(m)$ is chosen large enough, then
 \beq \label{eq11}
\big\| \mP^{\mC_{k(m)}} - \mP^\mC \big\|_{L^\infty (j\mC)} \ \leq \
\e.
 \eeq
Consider now the sequence $\big( \mF^{\mC_{k(m),l}}_{l}
\big)_{l=1}^\infty$ of Folding maps generated by the obtuse
simplices $\mC_{k(m),l}$ which approximate $\mC_{k(m)}$. The
sequence $\big( \mF^{\mC_{k(m),l}}_{l} \big)_{l=1}^\infty$ satisfies
estimates \eqref{eq7}, \eqref{eq8} and by compactness has a locally
uniform limit $\mF^{\mC_{k(m)}} \in [L^\infty_{loc}(\R^N)]^N$:
 \beq \label{eq13}
\mF^{\mC_{k(m),l}}_{l}\ \larrow \ \mF^{\mC_{k(m)}}
 \eeq
in $[L^\infty_{loc}(\R^N)]^N$ as $l \ri \infty$, along a
subsequence. We claim that $\mF^{\mC_{k(m)}} = \mP^{\mC_{k(m)}}$.
For, observe that since the boundary of $\mC_{k(m),l}$ consists of
tangent hyperplanes to the $C^1$ convex set $\mC_{k(m)}$ and $l$ is
the number of sides, there exists a dense sequence $(u_i)_1^\infty
\sub \p \mC_{k(m)}$, a respective sequence of outwards pointing
normal vectors $(\hat{u}_i)_1^\infty$ and a strictly increasing
function $\sigma \in C^0[0,\infty)$ with $\sigma(0)=0$ such that
 \beq \label{eq12}
\mF^{\mC_{k(m),l}}_l \big(u_i\ +\ t \hat{u}_i  \big)\ =  \ u_i\ + \
\sigma\big(2^{-l}\big) \hat{u}_i
 \eeq
for all $i\leq l$ and $t \in [0,\sigma(l)]$. Since the set of
half-lines
 \beq
D\ :=\ \bigcup_{i=1}^\infty\Big(u_i\ + \ \big\{ t\hat{u}_i : t\geq 0
\big\}\Big)
 \eeq
is dense in $\R^N \set \mC$, by passing to the limit to \eqref{eq12}
as $l\ri \infty$, we obtain
 \begin{align}
\mF^{\mC_{k(m)}} \big(u_i\ +\ t \hat{u}_i  \big)\ & =  \ u_i \\
 &= \ \mP^{\mC_{k(m)}} \big(u_i\ +\ t \hat{u}_i  \big), \nonumber
 \end{align}
for all $i\in \N$ and all $t \in [0,+\infty)$. Consequently,
$\mF^{\mC_{k(m)}} = \mP^{\mC_{k(m)}}$ on the dense set $D$; equality
on $\R^N\set\mC$ follows by continuity. Thus, by \eqref{eq13} we can
choose an $l(m) \in \N$ sufficiently large to get
 \beq \label{eq14}
\big\|\mF^{\mC_{k(m),l(m)}}_{l(m)} -  \mP^{\mC_{k(m)}}
\big\|_{L^\infty (j\mC)} \ \leq \ \e.
 \eeq
Finally, by \eqref{eq15} and \eqref{FIV}, we can increase $l(m)$
further to assure
 \beq \label{eq16}
\big\|\mF^{\mC_{k(m),l(m)}}_{l(m)} -  \mF^{\mC} \big\|_{L^\infty
(j\mC)} \ \leq \ \e.
 \eeq
By putting estimates \eqref{eq9}, \eqref{eq11}, \eqref{eq14} and
\eqref{eq16} together, the desired conclusion $\mF^{\mC} \equiv
\mP^{\mC}$ follows by the arbitrariness of $\e>0$ and $j\in \N$.
\qed

\ms

\noi The proof of Theorem \ref{Construction of the Approximations}
is now complete. \qed

\ms

\section{The Maximum Principles.}
\label{The Folding Principles.}

In this section we present our main results, the Maximum Principles
for approximate minimizers of nonconvex functionals, in the sense of
the Convex Hull Property. They are obtained by employing the
approximating Folding Maps of Theorem \ref{Construction of the
Approximations} to generate appropriate energy comparison functions
which compensate the lack of weak Lower Semi-Continuity. We shall
refer to them as \textit{Folding Principles}, since they are
obtained by employing the Folding Maps of Section \ref{Construction
of the Approximations.}.

Before proceeding to our main results, we first need to formulate
the notions of ``$u(\Om)$'' and ``$u(\p \Om)$'' in a weak sense,
meaningful for merely measurable functions. Let $u:\Om \sub \R^n
\larrow \R^N$ be a measurable function defined on the open set $\Om$
and let $A$, $K \sub \overline{\Om}$ be measurable as well. Let also
$|\_ |$ denote the Lebesgue measure in any dimension. If $|A|>0$,
the \emph{essential range} $u(A)$ is the closed set
 \beq \label{equ4}
u(A)\ := \ \left\{\eta \in \R^N\ \Big| \ \underset{x\in A}{\ess\,
\inf}\, |u(x)-\eta|=0\right\}.
 \eeq
If $|K|=0$, \eqref{equ4} is not directly meaningful for $K$. We
therefore define
 \beq \label{equ5}
u(K)\ := \ \bigcap\Big\{u(A)\ \Big| \ A\supseteq K, \ |A|>0\Big\}.
 \eeq
Finally, if $S\sub \R^N$, we denote the open $\e$-neighborhood of
$S$ by
 \beq
S^\e\ := \ \big\{\eta \in \R^N \ \big| \ |\eta - s|<\e ,\ s\in
S\big\}.
 \eeq

\subsection{The case of integrand with no direct dependence on the
function.} \label{subsec1}

\noi We consider first the simpler case of the functional
 \beq \label{equ6}
E(u,\Om)\ = \ \int_\Om \msL\big(x,Du(x)\big)dx.
 \eeq
\noi {\bf Hypotheses on the functional \eqref{equ6}.} We shall need
to impose the following assumptions: let $\Om \sub \R^n$ be open and
$\msL : \Om \by \R^{N \by n} \larrow \R$ a Carath\'eodory function
such that

\begin{enumerate}
\item \label{A} \emph{$\msL(x,\_ )$ is almost invariant under
the Orthogonal Group}: there exists $a\in L^1(\Om)$ such that for
a.e.\ $x \in \Om$, all $O \in O(N,\R)$ and all $P\in \R^{N\by n}$,
we have
 \[
\big|\msL(x,P)\ - \ \msL(x,OP)\big|\ \leq\ a(x).
 \]
\item \label{B} \emph{$\msL(x,\_ )$ is of $q$-growth}: there exist $C>0$,
$q\geq1$ and $b\in L^1(\Om)$ such that
 \[
- b(x) \ \leq \ \msL(x,P) \ \leq\ C\, |P|^q\ + \ b(x).
 \]
\end{enumerate}

\begin{remark} Assumption \eqref{B} is rather standard. Assumption
\eqref{A} is always satisfied if $\msL(x,\_ )$ is invariant under
$O(N,\R)$, in which case $a\equiv 0$. In particular, this always
holds in the frequent case where $\msL$ depends on $P$ through its
modulus: $\msL = \msL(x,|P|)$.

In the scalar case of $N=1$, assumption \eqref{A} requires
 \beq \label{equ9}
\big|\msL(x,P)\ - \ \msL(x,-P)\big|\ \leq\ a(x),
 \eeq
which means that $\msL(x,\_ )$ is almost even. Even in the scalar
case, evenness which occurs for $a\equiv 0$ and reads
$\msL(x,P)=\msL(x,-P)$, is an assumption much weaker than requiring
to be radial, i.e.\ $\msL = \msL(x,|P|)$.

We note that the form $\msL = \msL(x,|P|)$ of the integrand coupled
by monotonicity of the function $t \mapsto \msL(x,t)$ is a standing
assumption in the literature for the derivation of the Convex Hull
Property. In this work we weaken it substantially and in particular
we bypass quasiconvexity.

\end{remark}

\begin{theorem}[Folding Principle I] \label{FP1} Let $\Om \sub \R^n$ be open and
$\msL : \Om \by \R^{N \by n} \larrow \R$ a Carath\'edory function
satisfying assumptions \eqref{A} and \eqref{B}. Then, for any $\al>
0$ and $g \in [W^{1,q}(\Om) \cap L^\infty(\Om)]^N$ there exists an
$(\al+\|a\|_{L^1(\Om)})$-minimizer $u$ of the functional
\eqref{equ6} in $[W^{1,q}_g(\Om)]^N$ which satisfies
\begin{equation} \label{eq17}
u(\Om) \ \sub \ \co\, \big( u(\p \Om)^\al \big).
\end{equation}
\end{theorem}

\ms

By Theorem \ref{FP1} we readily obtain the following

\begin{corollary} In the setting of Theorem \ref{FP1}, if moreover
$\msL$ is invariant under $O(N,\R)$, that is $a\equiv 0$ in
assumption \eqref{A}, then there exists a minimizing family
$\{u_\al\}_{\al>0}$ of problem \eqref{equ3} which asymptotically
satisfies the Convex Hull Property as $\al \ri 0$:
\begin{equation} \label{eq18}
u_\al(\Om) \ \sub \ \co\, \big( u_\al(\p \Om)^\al \big).
\end{equation}
\end{corollary}

\begin{remark} \label{rem} By employing standard results (cf.\ for example Dacorogna
\cite{D1}, \cite{D3}), \eqref{eq18} is inherited by the
\emph{relaxed solution of \eqref{equ3}} and by \emph{minimizing
solutions of the Euler-Lagrange system of PDE}, if they exist.
\end{remark}

\begin{remark} In the absence of weak Lower Semi-Continuity,
\eqref{eq17} and \eqref{eq18} is \emph{all} we can assert, since we
can not pass to the limit as $\al \ri 0$ to sharpen them. Moreover,
if $\msL$ is not invariant under the Orthogonal Group, then the
Convex Hull Property is generally satisfied by an approximate
minimizer at a \emph{higher energy level}, increased by the amount
``$\|a\|_{L^1(\Om)}$'' of deviation of $\msL$ from invariance.
\end{remark}

\BPT \ref{FP1}. Fix $\al>0$. Under the lower $L^1$-bound of
assumption \eqref{B}, minimizing sequences of problem \eqref{equ3}
are equivalent to families of approximate minimizers as $\al \ri 0$.
Hence, we can choose an $\frac{\al}{2}$-minimizer $v \in
[W_g^{1,q}(\Om)]^N$ of \eqref{equ6}. We consider the
$\frac{\al}{2}$-neighborhood of $u(\p \Om)$ and set
 \beq \label{equ13}
\mC \ := \ \co\, \big( u(\p \Om)^{\al/2} \big).
 \eeq
Since $u-g \in [W^{1,q}_0(\Om)]^N$ and $g \in [L^{\infty}(\Om)]^N$,
by \eqref{equ4} and \eqref{equ5}, $\mC$ is a compact convex of
$\R^N$ with nonempty interior. By a translation, we may assume that
$0 \in \inter(\mC)$. Let $\mF_m$, $m\in \N$, be the Folding Maps of
Theorem \ref{Construction of the Approximations} relative to $\mC$.
Let also $v^m$ be the truncate of $v$ whose range is contained in
the largest ball centered at zero inside the $m$-dilate $m\, \mC$ of
$\mC$:
 \beq \label{equ10}
v^m \ := \ v\, \chi_{\{|v|\leq R(m)\}} \ + \ \frac{v}{|v|}\,
\chi_{\{|v|>R(m)\}},
 \eeq
 \beq \label{equ11}
R(m)\ := \ \max \Big\{R \in \N\ \Big| \ v\big(\{|v|\leq R\}\big) \,
\sub \, m\, \mC \Big\}.
 \eeq
Then, $v^m \ri v$ in $[W^{1,q}(\Om)]^N$ and $Dv^m \ri Dv$ a.e.\ on
$\Om$ as well, along a certain subsequence as $m\ri \infty$.
Moreover, by Theorem \ref{Construction of the Approximations}, the
transformation $\mF_m \circ v^m $ is well-defined and contained in
$[W^{1,q}(\Om)]^N$. Moreover, for $m=m(\al)$ large, in view of
\eqref{eq1} it satisfies
 \begin{align} \label{eq20}
\big(\mF_m \circ v^m\big)(\Om) \ & \sub \ \left(1+\frac{1}{m} \right)\mC \nonumber \\
 & \sub \ \mC^{\al/2}  \\
 & = \ \co\, \big( u(\p \Om)^{\al/2} \big)^{\al/2} \nonumber\\
 & = \ \co\, \big( u(\p \Om)^{\al} \big). \nonumber
 \end{align}
By \eqref{eq1}, $\mF_m$ leaves $\mC$ invariant; hence, by
\eqref{equ4}, \eqref{equ5} and \eqref{equ13} we obtain
 \beq \label{eq29}
v(\p \Om)\ = \ \big(\mF_m \circ v^m\big)(\p \Om).
 \eeq
Thus, $\mF_m \circ v^m \in [W^{1,q}_g(\Om)]^N$ and by \eqref{eq20}
and \eqref{eq29} we obtain
 \beq \label{eq26}
\big(\mF_m \circ v^m\big)(\Om) \ \sub \ \co \, \big(\mF_m \circ
v^m\big)(\p \Om).
 \eeq
By assumption \eqref{B} and since $|Dv^m|\leq 2|Dv|$ for $m$ large,
we have the bound
 \beq
\Big|\msL\big(x,Dv^m(x)\big)\Big|\ \leq \ 2^q C\, |Dv(x)|^q\ + \
b(x),
 \eeq
valid for a.e.\ $x\in \Om$. Hence, since $Dv^m \ri Dv$ a.e.\ as
$m\ri \infty$, the Dominated Convergence theorem implies
 \beq \label{eq23}
\int_\Om\msL\big(x,Dv^m(x)\big)dx\ \larrow \
\int_\Om\msL\big(x,Dv(x)\big)dx,
 \eeq
as $m\ri \infty$. By Theorem \ref{Construction of the
Approximations}, for a.e.\ $x\in \Om$, we have that
$D\mF_m\big(v^m(x)\big) \in O(N,\R)$. Hence, by employing assumption
\eqref{A} and \eqref{eq23}, we have
 \begin{align} \label{eq22}
E\big(\mF_m \circ v^m , \Om\big)\ &= \ \int_\Om \msL\big(x,D(\mF_m \circ v^m)(x)\big)dx  \nonumber\\
                &= \ \int_\Om \msL\big(x,D\mF_m(v^m(x)) \, Dv^m(x)\big)dx  \nonumber\\
                  &\leq \ \int_\Om \msL\big(x,Dv^m(x)\big)dx \ + \
                  \int_\Om a(x)\, dx\\
                  & \leq \  \int_\Om \msL\big(x,Dv(x)\big)dx\ + \frac{\al}{2} \ + \
                  \|a\|_{L^1(\Om)}  \nonumber\\
            & = \  E(v,\Om)\ + \frac{\al}{2} \ + \
                  \|a\|_{L^1(\Om)}, \nonumber
\end{align}
for $m=m(\al)$ large enough. Since $v$ is an
$\frac{\al}{2}$-minimizer of \eqref{equ6} and $v -\mF_m \circ v^m
\in [W^{1,q}_0(\Om)]^N$, by choosing $\psi \in [W^{1,q}_0(\Om)]^N$
arbitrary and setting
 \beq
\phi\ :=\ \psi - v + \mF_m \circ v^m,
 \eeq
we obtain that $\phi \in [W^{1,q}_0(\Om)]^N$. Thus, by approximate
minimality of $v$, we obtain
 \begin{align} \label{eq27}
E\big(v, \Om\big)\ &\leq \  E(v+\phi,\Om)\ + \frac{\al}{2}
\nonumber\\
&=\  E\big(\mF_m \circ v^m +\psi,\Om\big)\ + \frac{\al}{2}.
\end{align}
By combining \eqref{eq22} and \eqref{eq27}, we see that for any
$\al>0$, the function $u := \mF_{m} \circ  v^m$ is an $(\al
+\|a\|_{L^1(\Om)})$-minimizer which in view of \eqref{eq26}
satisfies \eqref{eq17}, as desired.
        \qed

\subsection{The general case of integrand which depends on
all the arguments.} \label{subsec1}

\noi Let us now consider the general functional
 \beq \label{equ29}
E(u,\Om)\ = \ \int_\Om \msL\big(x,u(x),Du(x)\big)dx.
 \eeq

\noi {\bf Hypotheses on the functional \eqref{equ29}.} We shall need
to impose the following assumptions: let $\Om \sub \R^n$ be open and
$\msL : \Om \by \R^N \by \R^{N \by n} \larrow \R$ a Carath\'eodory
function such that

\begin{enumerate}
\item \label{C} \emph{$\msL(x,\eta,\_ )$ is almost invariant under the
Orthogonal Group, locally uniformly in $\eta$}: for any $R>0$, there
exists an $a=a_R \in L^1(\Om)$ such that for a.e.\ $x \in \Om$, all
$|\eta|\leq R$, all $O \in O(N,\R)$ and all $P\in \R^{N\by n}$, we
have
 \[
\big|\msL(x,\eta,P)\ - \ \msL(x,\eta,OP)\big|\ \leq \ a(x).
 \]

\item \label{D} \emph{$\msL(x,\eta,\_ )$ is of $q$-growth}: there exist $C>0$,
$q\geq1$, $b\in L^1(\Om)$ and $d : \Om \by \R^N \larrow \R$ a
Carath\'eodory function such that
 \[
- b(x) \ \leq \ \msL(x,\eta,P) \ \leq\ C\, |P|^q\ + \ d(x,\eta).
 \]

\item \label{E} \emph{There exists a convex set such that the values of
$\msL(x,\_ ,P)$ outside of it almost exceed those on the boundary}:
there exist $l\in L^1(\Om)$ and $\mC\sub \R^N$ compact and convex
with $0 \in \inter(\mC)$ such that
 \[
 \max_{\eta \in \p \mC} \, \msL(x,\eta,P)\ \leq \
  \inf_{\eta \in \R^N \set \mC} \, \msL(x,\eta,P)\
 +\ l(x),
 \]
for a.e.\ $x\in \Om$ and all $P\in \R^{N\by n}$.
\end{enumerate}

\begin{remark} Assumptions \eqref{C} and \eqref{D} are analogous to
those of the previous case of \eqref{equ6}. Assumption \eqref{E}
says that there exists a convex set $\mC$ such that the function
$\msL(x,\_ ,P)$ has values on $\R^N \set \mC$ which \textit{almost}
exceed those on $\p \mC$. The nonnegative function $l$ represents
the deviation of the restriction $\msL(x,\_ ,P)\big|_{\mC}$ from
being constant. If $l\equiv0$, then ``almost'' in the previous
statement can be dropped. Assumption \eqref{E} is weaker than
requiring $\mC$ to be a sublevel set, since the values inside $\mC$
may well exceed those on the boundary.
\end{remark}

\begin{theorem}[Folding Principle II] \label{FP2} Let $\Om \sub \R^n$ be open and
$\msL : \Om \by \R^N \by \R^{N \by n} \larrow \R$ a Carath\'edory
function satisfying assumptions \eqref{C}, \eqref{D} and \eqref{E}.
Then, for any $\al> 0$ and any $g \in [W^{1,q}(\Om) \cap
L^\infty(\Om)]^N$ there exists an $(\al + \|a+l\|_{L^1(\Om)})
$-minimizer $u$ of the functional \eqref{equ29} in
$[W^{1,q}_g(\Om)]^N$ which satisfies
\begin{equation} \label{equ32}
u(\p \Om) \ \subset \ \mC \ \ \Longrightarrow \ \ u(\Om) \ \sub \
\mC^\al.
\end{equation}
\end{theorem}

By Theorem \ref{FP2} we readily obtain the following

\begin{corollary} In the setting of Theorem \ref{FP2}, if moreover
$a=l\equiv 0$ in assumptions \eqref{C} and \eqref{E}, then
there exists a minimizing family $\{u_\al\}_{\al>0}$ of problem
\eqref{equ3} which satisfies as $\al\ri 0$ that
\begin{equation} \label{equ33}
u_\al(\p \Om) \ \subset \ \mC \ \ \Longrightarrow \ \ u_\al(\Om) \
\sub \ \mC^\al.
\end{equation}
\end{corollary}

\begin{remark} A statement analogous to that or Remark \ref{rem} applies here as well.
Once again we observe that in general if $a\not\equiv0$ or
$l\not\equiv0$, there must be an increase at the energy level of the
resulting approximate minimizer.
\end{remark}

 \begin{remark}
The strictness assumption ``$u(\p \Om) \subset \mC$'' appearing in
\eqref{equ32} and \eqref{equ33} can be weakened to ``$u(\p \Om) \sub
\mC$'' by assuming in addition that \emph{the function $\eta \mapsto
\msL(x,\eta,P)$ has moduli of continuity over compacts depending
uniformly in $(x,P)\in \Om\by \R^{N\by n}$}. This is a rather weak
requirement which is always satisfied in the decoupled case of
$\msL(x,\eta,P) = A(x,P)+ W(\eta)$, a standard example of which is
given by the \emph{Action functional} $E(u,\Om)=\int_\Om
[\frac{1}{2}|Du(x)|^2+W(u(x))]dx$.
\end{remark}

\BPT \ref{FP2}. Fix $\al>0$ and let $v \in [W_g^{1,q}(\Om)]^N$ be an
$\frac{\al}{2}$-minimizer of \eqref{equ29} satisfying $v(\p \Om)
\subset \mC$. Let us denote the Folding Maps of $\mC$ provided by
Theorem \ref{Construction of the Approximations} by $\mF_m$, $m\in
\N$. Let $\mP$ denote the projection map onto $\mC$. Let finally
$v^m$ be the truncate of $v$ given by \eqref{equ10} and
\eqref{equ11}. Then, $\mF_m\circ v^m$ is a well defined function in
$[W^{1,q}(\Om)]^N$. By assumption \eqref{E}, $\mC$ is compact. By
\eqref{equ4} and \eqref{equ5}, $v(\p \Om)$ is closed and since $v(\p
\Om) \subset \mC$, there exists an open neighborhood of $v(\p \Om)$
which is contained into $\mC$. Hence, $\mF_m$ leaves the boundary
values invariant; we therefore obtain
 \beq
\big(\mF_m\circ v^m\big)(\p \Om) \ =\ v(\p \Om).
 \eeq
Consequently, we have that $\mF_m \circ v^m \in [W^{1,q}_g(\Om)]^N$.
Again by Theorem \ref{Construction of the Approximations}, for
$m=m(\al)$ large we have
 \begin{align} \label{eq43}
\big(\mF_m \circ v^m\big)(\Om)\ \sub \ \mC^\al.
 \end{align}
By \eqref{equ10}, \eqref{equ11} and Theorem \ref{Construction of the
Approximations}, we obtain that
 \begin{align}
Dv^m \ & \larrow\ Dv, \label{eq44a} \\
\mF_m \circ v^m \ &\larrow \ \mP \circ v, \label{eq44b}
 \end{align}
a.e.\ on $\Om$, both along a common subsequence as $m\ri \infty$. By
employing the bounds
 \begin{align} \label{eq45}
 \big|\mF_m \circ v^m\big|\ &\leq \ 2\,
\diam(\mC),\\
\big|\mP\circ v\big|\ &\leq \ \diam(\mC),  \label{eq46}\\
|Dv^m|\ & \leq \ 2\, |Dv|, \label{eq47}
 \end{align}
valid for $m$ large, by the bound of assumption \eqref{D} and
application of the Dominated Convergence theorem, we obtain
 \beq
\int_\Om \msL\Big(x, \mF_m \big( v^m(x)\big) ,Dv^m(x)\Big)dx \
\larrow \ \int_\Om \msL\Big(x, \mP \big(v(x)\big),Dv(x)\Big)dx,
\label{eq48}
 \eeq
as $m\ri \infty$. Now we employ the identity
 \begin{align} \label{eq50a}
E\big(\mF_m \circ v^m , \Om\big)\ &= \ \int_\Om \msL\Big(x,
                \big(\mF_m \circ v^m\big)(x),
                D\big(\mF_m \circ v^m\big)(x)\Big)dx  \nonumber\\
                &= \ \int_\Om \msL\Big(x, \mF_m \big( v^m(x)\big),
                                       D\mF_m \big(v^m(x)\big) \,
                                       Dv^m(x)\Big)dx
\end{align}
in order to estimate the energy $E\big(\mF_m \circ v^m , \Om\big)$.
For, by employing Theorem \ref{Construction of the Approximations},
for a.e.\ $x\in \Om$, we have that $D\mF_m\big(v^m(x)\big) \in
O(N,\R)$. We utilize \eqref{eq48} and assumption \eqref{C} where as
$R$ we take
 \beq
R\ :=\ 2\, \diam(\mC)
 \eeq
with respective $a=a_R$. Then, identity \eqref{eq50a} implies
 \begin{align} \label{eq50b}
E\big(\mF_m \circ v^m , \Om\big)\  &\leq \ \int_\Om \msL\Big(x,
\mF_m \big(
                  v^m(x)\big),Dv^m(x)\Big)dx \nonumber\\
                  &\ \ \ \ + \ \int_\Om a(x)\, dx  \nonumber\\
                  & \leq \  \int_\Om \msL\Big(x,\mP(v(x)),Dv(x)\Big)dx\
                  + \frac{\al}{2} \\
                  &\ \ \ \ + \ \|a\|_{L^1(\Om)},  \nonumber
\end{align}
for $m=m(\al)$ large enough. By assumption \eqref{E}, a.e.\ on the
set
 \beq
\{v \not\in \mC \}\ := \ \Big\{ x\in \Om\ \Big|\ v(x)\in \R^N \set
\mC \Big\} \ \sub \ \Om
 \eeq
we have the estimate
 \begin{align} \label{eq51}
 \msL\big(\_ ,\mP \circ v,Dv\big) \ & \leq \ \max_{\eta \in \p \mC}
                                            \msL\big(\_ ,\eta,Dv\big) \nonumber\\
                                    & \leq \ \inf_{\eta \in \R^N \set \mC}
                                           \msL\big(\_ ,\eta,Dv\big)   \ + \ l \\
                                    & \leq \ \msL\big(\_ ,v,Dv\big)  \ + \
                                    l.\nonumber
 \end{align}
Trivially, a.e.\ on its complement $\{v \in \mC \}$ we have
 \beq  \label{eq52}
 \msL\big(\_ ,\mP \circ v,Dv\big)  \ = \  \msL\big(\_ ,v,Dv\big).
 \eeq
By employing \eqref{eq51} and \eqref{eq52}, estimate \eqref{eq50b}
implies
 \begin{align} \label{eq53}
E\big(\mF_m \circ v^m , \Om\big)\ & \leq \  \int_{\{v \not\in \mC
                                   \}} \msL\Big(x,\mP(v(x)),Dv(x)\Big)dx \nonumber\\
                                    &\ \ \ \ + \ \int_{\{v \in \mC \}}
                                   \msL\Big(x,\mP(v(x)),Dv(x)\Big)dx \nonumber\\
                                &\  \ \ \ + \frac{\al}{2} \ + \ \|a\|_{L^1(\Om)}
                                   \\
                                    & \leq \  \int_{\{v \not\in \mC
                                   \}} \msL\Big(x,v(x),Dv(x)\Big)dx\
                                   + \ \int_\Om l(x)dx \nonumber\\
                                    &\ \ \ \ + \ \int_{\{v \in \mC \}}
                                   \msL\Big(x,v(x),Dv(x)\Big)dx \nonumber\\
                                &\  \ \ \ + \frac{\al}{2} \ + \ \|a\|_{L^1(\Om)}
                                   \nonumber.
\end{align}
Hence, by \eqref{eq53} we obtain
\begin{align} \label{eq53a}
E\big(\mF_m \circ v^m , \Om\big)\ & \leq \  \int_\Om
\msL\Big(x,v(x),Dv(x)\Big)dx \nonumber\\
&\ \ \ \  + \frac{\al}{2} \ + \
\|a\|_{L^1(\Om)}\ + \ \|l\|_{L^1(\Om)} \\
                 & = \ E(v,\Om) \ + \frac{\al}{2} \ + \
\|a+l\|_{L^1(\Om)}. \nonumber
\end{align}
Consequently, since $v$ is an $\frac{\al}{2}$-minimizer of
\eqref{equ29} and $v -\mF_m \circ v^m \in [W^{1,q}_0(\Om)]^N$, by
choosing $\psi \in [W^{1,q}_0(\Om)]^N$ arbitrary and setting
 \beq
\phi\ :=\ \psi - v + \mF_m \circ v^m,
 \eeq
we obtain that $\phi \in [W^{1,q}_0(\Om)]^N$. Thus, by approximate
minimality of $v$, we obtain
 \begin{align} \label{eq54}
E\big(v, \Om\big)\ &\leq \  E(v+\phi,\Om)\ + \ \frac{\al}{2}
\nonumber\\
&=\  E\big(\mF_m \circ v^m +\psi,\Om\big)\ +\ \frac{\al}{2}.
\end{align}
Hence, by combining \eqref{eq53a}, \eqref{eq54} and \eqref{eq43},
for any $\al>0$ the function $u := \mF_{m} \circ v^m$ is an $(\al
+\|a+l\|_{L^1(\Om)})$-minimizer satisfying property \eqref{equ32}.
Theorem \ref{FP2} follows.
  \qed

\ms

\medskip
\noindent {\bf Acknowledgement.} {This work was done when the Author was a doctoral student at the Department of Mathematics, University of Athens, Greece. 

I wish to thank  L.\ Ambrosio, J. Ball, P.\ Bates, B.\ Dacorogna, L.C.\ Evans, F. Leonetti, and A.\ Tertikas
for their suggestions and comments on an earlier weaker version of
this manuscript which led to substantial improvements.

Last but not least, I am profoundly indebted to N.\ Alikakos for his
constructive criticism and for our inspiring scientific
discussions.}

\ms

\footnotesize 

\end{document}